\documentclass[12pt]{amsart}
\usepackage{amsfonts,amscd,amsmath,amssymb,amsthm,mathrsfs}
\usepackage{array,latexsym,float,enumerate}
\usepackage{graphicx,epsfig}
\usepackage[top=3cm, bottom=3cm, left=2.5cm, right=2.5cm,twoside=false]{geometry}
\usepackage[all]{xy}
\numberwithin{equation}{section}

\theoremstyle{plain}
\newtheorem{theorem}{Theorem}[section]
\newtheorem{lemma}[theorem]{Lemma}
\newtheorem{proposition}[theorem]{Proposition}
\newtheorem{corrolary}[theorem]{Corollary}
\theoremstyle{definition}
\newtheorem{definition}[theorem]{Definition}
\newtheorem{example}[theorem]{Example}
\newtheorem{noname}[theorem]{}
\newtheorem{remark}[theorem]{Remark}
\newtheorem{construction}[theorem]{Construction}
\newtheorem{notation}[theorem]{Notation}
\theoremstyle{remark}
\newtheorem*{smallremark}{Remark}

\newtheorem{case}{Case} \makeatletter \@addtoreset{case}{theorem}\makeatother
\newtheorem{claim}{Claim} \makeatletter \@addtoreset{claim}{theorem}\makeatother

\newcommand{\bthm}{\begin{theorem}}
\newcommand{\bprop}{\begin{proposition}}
\newcommand{\blem}{\begin{lemma}}
\newcommand{\bcor}{\begin{corrolary}}
\newcommand{\brem}{\begin{remark}}
\newcommand{\bdfn}{\begin{definition}}
\newcommand{\bitem}{\begin{itemize}}
\newcommand{\bex}{\begin{example}}
\newcommand{\bno}{\begin{noname}}
\newcommand{\bsrem}{\begin{smallremark}}
\newcommand{\bnot}{\begin{notation}}
\newcommand{\bcon}{\begin{construction}}
\newcommand{\bca}{\begin{case}}
\newcommand{\bcl}{\begin{claim}}

\newcommand{\ecl}{\end{claim}}
\newcommand{\eca}{\end{case}}
\newcommand{\econ}{\end{construction}}
\newcommand{\enot}{\end{notation}}
\newcommand{\esrem}{\end{smallremark}}
\newcommand{\eno}{\end{noname}}
\newcommand{\eex}{\end{example}}
\newcommand{\eitem}{\end{itemize}}
\newcommand{\ethm}{\end{theorem}}
\newcommand{\eprop}{\end{proposition}}
\newcommand{\elem}{\end{lemma}}
\newcommand{\ecor}{\end{corrolary}}
\newcommand{\erem}{\end{remark}}
\newcommand{\edfn}{\end{definition}}
\newcommand{\benum}{\begin{enumerate}}
\newcommand{\eenum}{\end{enumerate}}

\newcommand{\wt}{\widetilde}
\newcommand{\cal}[1]{\mathcal{#1}}

\newcommand{\ds}{\displaystyle}

\def\8{\infty}
\def\.{\cdot}
\def\PP{\mathbb{P}}

\def\Q{\mathbb{Q}}

\def\:{\colon}

\def\bsk{\bigskip}

\def\Bk{\operatorname{Bk}}

\def\Supp{\operatorname{Supp}}
\def\Pic{\operatorname{Pic}}

\usepackage[hypertex,linkbordercolor={0 0 1}]{hyperref} 

\begin{document}

\title[Coolidge-Nagata conjecture]{The Coolidge-Nagata conjecture holds \\ for curves with more than four cusps}

\author[Karol Palka]{Karol Palka}
\address{Karol Palka: D\'epartement de math\'{e}matiques, Universit\'{e} du Qu\'{e}bec \`{a} Montr\'{e}al and Institute of Mathematics, Polish Academy of Sciences, Warsaw}

\email{palka@impan.pl}

\begin{abstract} Let $E\subseteq \PP^2$ be a rational curve defined over complex numbers which has only locally irreducible singularities. The Coolidge-Nagata conjecture states that $E$ is \emph{rectifiable}, i.e. it can be transformed into a line by a birational automorphism of $\PP^2$. We show that if it is not rectifiable then the tree of the exceptional divisor for its minimal embedded resolution of singularities has at most nine maximal twigs. This settles the conjecture in case $E$ has more than four singular points. \end{abstract}

\maketitle

\section{Main result}\label{sec:main result}

All varieties considered are complex algebraic. An irreducible curve is \emph{cuspidal} if and only if all its singular points are cusps, i.e. they are locally irreducible. We are interested in rational cuspidal curves embedded into the projective plane $\PP^2$. Let $\bar E\subseteq \PP^2$ be such a curve. There are many examples with $\bar E$ having one, two or three cusps (already among quartics). Up to a choice of coordinates on the plane there is only one known example with four cusps. It is of degree five and has parametrization $(t^3-1,t^5+2t^2,t)$. No examples with more than four cusps are known and it is expected that they do not exist. It is also expected that any rational plane cuspidal curve is \emph{rectifiable}, i.e. there exists a birational automorphism of the plane, such that the proper transform of the curve is a line. This is known as the \emph{Coolidge-Nagata conjecture/problem}. The conjecture has been verified for all known examples. We show that even if rational plane cuspidal curves with more than four cusps do exist, the Coolidge-Nagata conjecture necessarily holds for them.

\bthm\label{thm:CN_holds_for_c>4} Let $\bar E\subseteq \PP^2$ be a rational cuspidal curve defined over complex numbers. If $\bar E$ has more than four cusps then there exists a birational automorphism of $\PP^2$ which transforms $\bar E$ into a line. \ethm

We show if fact that if $\bar E\subseteq \PP^2$ is non-rectifiable then the tree of the exceptional divisor for its minimal embedded resolution has at most nine maximal twigs (cf. \ref{cor:nonrectifiable implies t<=9}).

The structure of the paper is as follows. In section \ref{sec:two inequalities} we prove two main inequalities. In section \ref{sec:5 cusps} we complete the proof of the theorem by dealing with the case of five cusps and ten maximal twigs. In section \ref{sec:four cusps} we exclude the case of four cusps and ten maximal twigs.

\bsk\textsl{\textsf{Acknowledgements.}} The author would like to thank Prof. Mariusz Koras for teaching him the theory of open surfaces and to Prof. Peter Russell and Prof. Steven Lu for inviting him to Montreal and providing great working conditions.

\tableofcontents

\section{Preliminaries}\label{sec:preliminaries}

\subsection{Open surfaces}\label{ssec:open surfaces} We recall some results from the theory of non-complete surfaces, also to settle the notation. For a complete treatment the reader is referred to \cite{Miyan-OpenSurf}.

Let $T$ be a nonzero reduced simple normal crossing divisor on a smooth complete surface $X$. We denote the Iitaka-Kodaira dimension of $T$ by $\kappa(T)$. If $R$ is a reduced divisor with support contained in $T$ we define $\beta_T(R)=R\cdot (T-R)$ and call it a \emph{branching number of $R$ in $T$}. If $R$ is irreducible we say that $R$ is a \emph{tip} or a \emph{branching component} if $\beta_T(R)=1$ or $\beta_T(R)\geq 3$ respectively.

The \emph{arithmetic genus of $T$} is $p_a(T)=\frac{1}{2}T\cdot (K+T),$ where $K$ is the canonical divisor (class) on $X$. $T$ is a \emph{rational tree} if all its components are rational and the dual graph of $T$ contains no loops. In this case $p_a(T)=0$. We call $T$ a \emph{chain} if it has no branching components. If $T=T_1+\ldots+T_k$, is a decomposition of a rational chain into irreducible components, such that $T_i\cdot T_{i+1}=1$ for $i<k$, then we write $T=[-T_1^2,-T_2^2,\ldots,-T_k^2].$ By $(m)_p$ we mean a sequence $(m,m,\ldots,m)$ of length $p$. An $(n)$-curve is a smooth rational curve with self-intersection $n$.

We define the \emph{discriminant of $T$} by $d(T)=\det(-Q(T))$, where $Q(T)=(T_i\cdot T_j)_{i,j\leq k}$ is the intersection matrix of $T$. We put $d(0)=1$. The following formula follows from elementary properties of determinants (cf. \cite{Russell2}).

\blem\label{lem:det formula} Let $S$ and $T$ be reduced simple normal crossing divisors, such that $S\cdot T=1$ and let $S_0\subseteq S$ and $T_0\subseteq T$ be the irreducible components for which $S_0\cdot T_0=1$. Then $$d(S+T)=d(S)d(T)-d(S-S_0)d(T-T_0).$$ In particular, if $S=S_0$ then $d(T+S_0)=-S_0^2d(T)-d(T-T_0)$. \elem

\blem\label{lem:chains with small d} Let $T$ be a rational chain which contains no $(-1)$-curves and has a negative definite intersection matrix. \benum[(i)]

\item If $d(T)=2$ then $T=[2]$.

\item If $d(T)=3$ then $T=[2,2]$ or $T=[3]$.

\item If $d(T)=4$ then $T=[2,2,2]$ or $T=[4]$.

\item If $d(T)=5$ ten $T=[2,2,2,2]$ or $T=[2,3]$ or $T=[3,2]$ or $T=[5]$.

\item If $d(T)=6$ then $T=[2,2,2,2,2]$ or $T=[6]$. \eenum\elem

\begin{proof} We just note that if $T_1$ is the tip of $T$ and $T_2$ is the component meeting $T_1$ then by \ref{lem:det formula} $d(T)=-T_1^2d(T-T_1)-d(T-T_1-T_2)>d(T-T_1)$, so all chains $T$ with given discriminant can be found by induction. \end{proof}

Assume now that $T$ is a rational tree without non-branching $(-1)$-curves and with intersection matrix which is not negative definite. Assume also that $T$ is not a chain and that the intersection matrices of all its maximal twigs are negative definite. Let $T_i=T_{i,1}+\ldots+T_{i,k_i},$ where $T_{i,1}$ assumed to be a tip of $T$, $i=1,\ldots,t$, be all its maximal twigs. We put $e(T_i)=d(T_i-T_{i,1})/d(T_i)$ and $\delta(T_i)=1/d(T_i)$. We define $$\delta(T)=\sum_{i=1}^s\delta(T_i) \text{\ \ and \ } e(T)=\sum_{i=1}^se(T_i).$$

Assume that $T$ is as above and $\kappa(K+T)\geq 0$. We have the Zariski decomposition $K+T=(K+T)^++(K+T)^-,$ where $(K+T)^+$ is numerically effective and $(K+T)^-$ is effective, either empty or having a negative definite intersection matrix. Moreover, $(K+T)^+\cdot B=0$ for any curve $B$ contained in $\Supp (K+T)^-$. We define $\Bk T$, the \emph{bark of $T$}, as a unique $\Q$-divisor with support contained in the sum of maximal twigs of $T$ and satisfying $$\Bk T\cdot T_0=\beta_T(T_0)-2$$ for every component $T_0$ of every maximal twig of $T$.

\blem\label{lem:Bk} Let $T$ be a rational tree as above. Let $T_i$ be a maximal twig of $T$ and let $T_0$ be a component of $T_i$. Denote the coefficient of $T_0$ in the decomposition of $\Bk T$ into irreducible components by $t_0$. \benum[(i)]

\item The coefficient $t_0$ satisfies $0<t_0<1$.

\item If $T_0$ meets $T-T_i$ then $t_0=\delta(T_i)$.

\item $(\Bk T)^2=-e(T)$.

\item If there is no $(-1)$-curve $A$ on $X$, for which $T\cdot A\leq 1$, then $(K+T)^-=\Bk T$.

\eenum \elem

\begin{proof} Write $T_i=T_{i,1}+T_{i,2}+\ldots+T_{i,k_i}$, where $T_{i,k}$ are irreducible and $T_{i,j}\cdot T_{i,j+1}=1$ for $j<k_i$. Then by \cite[2.3.3.4]{Miyan-OpenSurf} the coefficient of $T_{i,j}$ in $\Bk T$ equals $d(T_{i,j+1}+\ldots+T_{i,k_i})/d(T_i)$. This gives (i), (ii) and (iii). Part (iv) follows from 2.3.11 loc. cit. \end{proof}

\subsection{Rectifiability}\label{ssec:rectifiability} Let $\bar E\subseteq \PP^2$ be a rational curve with singular points $q_1,\ldots,q_c$, $c>0$. Let $\pi\:X\to \PP^2$ be a minimal embedded good resolution of singularities for $\bar E$, i.e. $\pi^*\bar E$ is a simple normal crossing divisor, whose all $(-1)$-curves are branching. Let $K$ be the canonical divisor (class) on $X$ and let $D$ be the reduced total transform of $\bar E$. We denote the proper transform of $\bar E$ on $X$ by $E$. Write $\pi^{-1}(q_i)=Q_i$, where $Q_i$ is reduced effective. By definition $Q_i$ is a rational tree with negative definite intersection matrix. Since $Q_i$ contracts to a smooth point on a surface, $d(Q_i)=d([1])=1$. Note that if $\bar E$ is cuspidal then $D$ is a rational tree.

\bprop\label{prop:known_results} Let $\bar E$ and $E$ be as above. \benum[(i)]

\item $\bar E\subseteq \PP^2$ is rectifiable if and only if $\kappa(K+E)=-\8$.

\item If $\kappa(K+E)\geq 0$ then $h^0(2K+E)\geq 1$. The inequality is strict if $\kappa(K+E)=2$.

\eenum \eprop

\begin{proof} (i) is a theorem of Coolidge \cite{Coolidge}, see also \cite[2.6]{Kumar-Murthy}. (ii), crucial for us, is \cite[2.4, 3.2]{Kumar-Murthy}. Let us recall the proof of the first part of (ii), rearranging the arguments a bit. First of all, we can assume without loss of generality that there is no $(-1)$-curve on $X$ disjoint from $E$. Let $m>0$ be minimal such that $h^0(mE+m'K)\neq 0$ for some $m'\geq m$. Write $$mE+m'K=\sum_i U_i,$$ where $U_i$'s are irreducible. By the choice of $m$, $U_i\neq E$. Note that if $(K+E)\cdot U_{i_0}<0$ for some $i_0$ then $K\cdot U_{i_0}<0$, so $$(\sum_i U_i)\cdot U_{i_0}=m(K+E)\cdot U_{i_0}+(m'-m)K\cdot U_{i_0}<0.$$ But in this case we would get $U_{i_0}^2<0$, so $U_{i_0}$ would be a $(-1)$-curve disjoint from $E$, which contradicts the assumption. We get $$0\leq (K+E)\cdot (mE+m'K)=-2m+m'K\cdot(K+E),$$ so $K\cdot (K+E)>0$. Since the numerical class of $K+E$ is nonzero, by the Riemann-Roch theorem $h^0(2K+E)\geq K\cdot (K+E)>0$. \end{proof}

\brem It follows from \ref{prop:known_results}(ii) that if $\kappa(K+E)\geq 0$ then $E^2\leq -4$ and $\deg \bar E\geq 6$ (and the inequality is strict if $\kappa(K+E)=2$). The see the former note that since $\kappa(K)=-\8$, $E$ is not in the fixed part of $|2K+E|$, so $0\leq E\cdot (2K+E)=-4-E^2$. The latter follows from the equality $\pi_*(2K+E)=2K_{\PP^2}+\bar E$.

Note also that if $\bar E\subseteq \PP^2$ is a general rational curve of degree $d$ then its singularities are ordinary double points (nodes), so we compute easily $2K+E\sim (d-6)H$, where $H$ is a pullback of a line on $\PP^2$. Thus $\bar E$ is not rectifiable for $d\geq 6$. A general rational sextic with ten nodes is an example of lowest degree. \erem

\blem\label{lem:how k(K+D) depends on c} Let $\bar E\subseteq \PP^2$ be rational and cuspidal. \benum[(i)]

\item If $c\geq 2$ then $h^0(2(K+D))\neq 0$.

\item If $c\geq 3$ then $\kappa(K+D)=2$.

\item If $\kappa(K+D)=2$ then $(K+D)^-=\Bk D$. \eenum \elem

Part (ii) is explicitly stated in \cite{Wakabayashi} and part (i) follows from a proof there.  For (iii) note that if $\kappa(K+D)=2$ then by the Lefschetz duality $X-D$ is a smooth $\Q$-acyclic surface of general type, so it does not contain topologically contractible lines by \cite{MiTs-lines_on_qhp}. Since $X-D$ is affine, it follows that there is no $(-1)$-curve $A$ on $X$ for which $A\cdot D\leq 1$, so $(K+D)^-=\Bk D$ by \ref{lem:Bk}(iv).

\section{Two inequalities} \label{sec:two inequalities}

From now on we assume that $\bar E\subseteq \PP^2$ is a rational cuspidal curve with cusps $q_1,\ldots,q_c$, $c>0$. In this case the divisor $Q_i$ can be seen as produced by a \emph{connected sequence of blow-ups}, i.e. we can decompose the morphism contracting $Q_i$ to a point into a sequence of blow-ups $\sigma_1\circ\ldots\circ\sigma_s$, so that then the center of $\sigma_{i+1}$ belongs to the exceptional component of $\sigma_i$ for $i\geq 1$. Let $C_i$ be the unique $(-1)$-curve in $Q_i$. Clearly, $E\cdot Q_i=E\cdot C_i=1$. Since $\pi$ is minimal, $C_i$ is not a tip of $D$, so $Q_i-C_i$ has two connected components. One of them is a rational chain and the other a rational tree, $C$ meets them in tips. We denote the maximal twigs of $D$ by $T_1,\ldots,T_t$.

The surface $X-D$ is $\Q$-acyclic, i.e. $H_i(X-D,\Q)=0$ for $i>0$. The group $(\Pic \PP^2)\otimes \Q$ is generated by $\bar E$, so $(\Pic X)\otimes \Q$ is generated freely by the components of $D$. Since $D-E$ has a negative definite intersection matrix, the intersection matrix of $D$ is not negative definite by the Hodge index theorem. If $\kappa(K+D)\geq 0$ we put $\cal P=(K+D)^+$. In the inequalities below the case $c=1$ is somewhat special. It is convenient to introduce $\epsilon(c)$ defined to be $0$ for $c>1$ and $1$ for $c=1$.

\bprop If $\kappa(K+E)\geq 0$ then the following inequality holds: \begin{equation*} t-\frac{1}{2}(c+\epsilon(c))\leq \cal \delta(D)+1+\cal P^2\leq \delta(D)+4.\tag{$\star$}\label{eq:*}\end{equation*}
\eprop

\begin{proof} We have $\kappa(K+D)\geq\kappa(K+E)\geq 0$. The divisor $R=D-T_1-\ldots-T_t$ is a reduced rational tree, so $p_a(R)=0$. The rational twigs are contained in $\Supp \Bk D$, so they intersect $\cal P$ trivially. We get $$\cal P\cdot D=\cal P\cdot R=(K+D-\Bk D)\cdot R= (K+R)\cdot R+ \beta_D(R)-\Bk D\cdot R.$$ By \ref{lem:Bk}(ii) the coefficient in $\Bk D$ of the component of the twig $T_i$ which intersects $R$ is $\delta(T_i)$. Thus $\cal P\cdot D=-2+t-\delta(D)$. We have $$\cal P\cdot E=(K+D-\Bk D)\cdot E=(K+E)\cdot E+\beta_D(E)-\Bk D\cdot E=-2+c-\Bk D\cdot E.$$ If $c=1$ then $E$ is a maximal twig of $D$, which gives $\cal P\cdot E=c-1=0$. If $c\geq 2$ then $E$ is disjoint from $\Bk D$, hence $\cal P\cdot E=c-2$. The formula which holds for both cases is therefore $\cal P\cdot E=c+\epsilon(c)-2$. From \ref{prop:known_results}(ii) we see that $$0\leq \cal P\cdot (2K+E)=2\cal P(K+D)-2\cal P\cdot D+
\cal P\cdot E=2\cal P^2-2\cal P\cdot D+\cal P\cdot E.$$ Thus $$\cal P^2\geq \cal P\cdot D-\frac{1}{2}\cal P\cdot E=t-2-\delta(D)-\frac{1}{2}(c+\epsilon(c)-2)=t-1-\delta(D)-\frac{1}{2}(c+\epsilon(c)).$$  Finally, from the logarithmic Bogomolov-Miyaoka-Yau inequality (see \cite{Langer}, cf. \cite[2.5]{Palka-exceptional}) we have $\cal P^2\leq 3\chi(X-D)=3$, where $\chi$ denotes the Euler characteristic.
\end{proof}

\bprop If $\kappa(K+D)=2$ then \begin{equation} h^0(2K+D)+e(D)=\cal P^2+2\leq 5. \tag{$\diamond$}\label{eq:diamond}\end{equation} \eprop

\begin{proof} We have $K\cdot (K+D)-2=(K+D)^2=\cal P^2+\Bk^2 D=\cal P^2 -e(D).$ Since the numerical class of $K+D$ is nonzero, the Riemann-Roch theorem gives $h^0(2K+D)-h^1(2K+D)=K\cdot (K+D)$. By the logarithmic Bogomolov-Miyaoka-Yau inequality we have $\cal P^2\leq 3$, so we see that $$h^0(2K+D)-h^1(2K+D)+e(D)=\cal P^2+2\leq 5.$$ Note that by \ref{lem:Bk}(i) $\Supp \Bk$ is an effective $\Q$-divisor with simple normal crossing support and proper fractional coefficients. Since $\cal P$ is nef and big ($\kappa(\cal P)=\kappa(K+D)=2$), the Kawamata-Viehweg vanishing theorem (see for example \cite[9.1.18]{Lazarsfeld_II}) says that $h^1(2K+D)=0$ for $i>0$. This completes the proof. \end{proof}

\bsrem Note that every $Q_i$ contains some maximal twig $T_j$ of $D$ with $d(T_j)\geq 3$. Indeed, it is clear if $Q_i=[(2)_m,3,1,2]$ for some $m$ and the general case follows by induction on the number of components of $Q_i$. \esrem

\bcor\label{cor:nonrectifiable implies t<=10} Assume $\bar E\subseteq \PP^2$ is not rectifiable. Then $t\leq 10$ in case $c\geq 4$ and $t\leq 9$ in case $c\leq 3$. In particular, $\bar E$ has at most five cusps. \ecor

\begin{proof} The remark above gives $\delta(D)\leq \frac{1}{3}c+ \frac{1}{2}(t-c)= \frac{1}{2}t-\frac{1}{6}c$. By \ref{prop:known_results}(i) $\kappa(K+E)\geq0$, so \eqref{eq:*} gives $t\leq \frac{2}{3}c+\epsilon(c)+8$. For $c\leq 2$ we get $t\leq 9$. Assume $c\geq 3$. By \ref{lem:how k(K+D) depends on c}(ii) $\kappa(K+D)=2$. Since $2c\leq t$, by \eqref{eq:*} and \eqref{eq:diamond} $$\frac{3}{4}t\leq t-\frac{1}{2}c\leq 4+\delta(D)\leq 4+e(D)\leq 8,$$ so $t\leq \frac{32}{3}<11$. If $c=3$ then we have in fact $t-\frac{3}{2}\leq 8$, so $t\leq 9$. \end{proof}

\bsrem We note that if $\bar E\subseteq \PP^2$ is rectifiable then it has at most eight cusps by \cite{Tono-on_the_number_of_cusps}. \esrem

\section{Five cusps}\label{sec:5 cusps}

In this section we assume that $\bar E\subseteq \PP^2$ is a rational cuspidal curve which is non-rectifiable. In particular $h^0(2K+D)\geq h^0(2K+E)>0$. By \ref{cor:nonrectifiable implies t<=10} $t\leq 10$, so $\bar E$ has at most five cusps. Therefore, to prove the theorem \ref{thm:CN_holds_for_c>4} we can assume that $t=10$ and $c=5$. For an (ordered) rational chain with negative definite intersection matrix we put $u(T)=e(T)-\delta(T)\geq 0$.

\blem\label{lem:properties of u} Assume $T=T'+C+T''$ is a rational chain with a negative definite intersection matrix, with a unique $(-1)$-curve $C$ and having $d(T)=1$. Assume $T',T''$ are nonempty and $d(T')\leq d(T'')$. Put $\bar u(T)=u(T')+u(T'')$ (we assume that the tips of $T$ are the first components of $T'$ and $T''$). Then: \benum[(i)]

\item $d(T')$ and $d(T'')$ are coprime,

\item $\bar u(T)\geq 0$ and the equality holds only if $T=[2,1,3]$,

\item if $T'=[2]$ then $T''=[(2)_k,3]$ for some $k\geq 0$ and $\bar u(T)=1-\frac{3}{2k+3}$,

\item if $T'=[2,2]$ then $T''=[(2)_k,4]$ for some $k\geq 0$ and $\bar u(T)=\frac{4}{3}-\frac{4}{3k+4}$,

\item if $T'=[3]$ then $T''=[(2)_k,3,2]$ for some $k\geq 0$ and $\bar u(T)=1-\frac{4}{3k+5}$,

\item if $T'=[4]$ then $T''=[(2)_k,3,2,2]$ for some $k\geq 0$ and $\bar u(T)=1-\frac{5}{4k+7}$.

\eenum \elem

\begin{proof} Let $A$ and $A'$ be the components of $T'$ and $T''$ meeting $C$ respectively. By \ref{lem:det formula} $d(T')d(C+T'')-d(T'-A')d(T'')=d(T)=1$, so $d(T')$ and $d(T'')$ are coprime. If $\bar u(T)=0$ then $T'$ and $T''$ are irreducible, so $d([-(T')^2,1,-(T'')^2])=1$, which happens only if $T=[2,1,3]$. For (iii) note that since $T$ contracts to a $(-1)$-curve, we necessarily have $T=[2,1,3,(2)_k]$ for some $k\geq 0$. We compute $d([(2)_k,3])=2k+3$, which gives the result. The remaining cases are done analogously. \end{proof}

\bcor\label{cor: chains with small u} Let $T$ be as above. Put $\bar \delta(T)=\delta(T')+\delta(T'')$. If $d(T')\leq 4$ and $0<\bar u(T)<\frac{1}{2}$ then:\benum[(i)]

\item $T=[2,1,3,2]$ and $(\bar u(T),\bar \delta(T))=(\frac{2}{5},\frac{7}{10})$ or

\item $T=[3,1,2,3]$ and $(\bar u(T),\bar \delta(T))=(\frac{1}{5},\frac{8}{15})$ or

\item $T=[2,2,1,4]$ and $(\bar u(T),\bar \delta(T))=(\frac{1}{3},\frac{7}{12})$ or

\item $T=[4,1,2,2,3]$ and $(\bar u(T),\bar \delta(T))=(\frac{2}{7},\frac{11}{28})$.  \eenum\ecor

\begin{proof} If $d(T')\leq 4$ then $T'=[2], [3], [4], [2,2]$ or $[2,2,2]$. In the last case $\bar u(T)\geq u([2,2,2])=\frac{1}{2}$. For the remaining cases use \ref{lem:properties of u}. \end{proof}

Put $\bar u_i=\bar u(Q_i)$ and $\bar \delta_i=\bar \delta(Q_i)$. We have $\delta(D)\neq e(D)$, otherwise  by \ref{lem:properties of u}(ii) $e(D)=5(\frac{1}{2}+\frac{1}{3})>4$, which contradicts \eqref{eq:diamond}. The inequalities \eqref{eq:*} and \eqref{eq:diamond} give \begin{equation} \frac{7}{2}\leq\delta(D)< e(D)\leq 4.\label{eq:t=10, bounds on delta and e}\end{equation} In particular $$0<\sum_{i=1}^{5} \bar u_i\leq \frac{1}{2}.$$ By renumbering twigs we can assume that $Q_i-C_i=T_i+T_{i+5}$ and $d(T_i)< d(T_{i+5})$.

\bprop\label{prop:only u5 is nonzero} With the notation as above $\bar u_1=\bar u_2=\bar u_3=0$ and: \benum[(i)]

\item either $Q_4=[3,1,2,3]$ and $Q_5=[3,1,2,3]$ or

\item $\bar u_4=0$ and $Q_5=[n,1,(2)_{n-2},3]$ for some $5\leq n\leq 9$ or

\item $\bar u_4=0$ and $Q_5=[5,2,1,3,2,2,3]$. \eenum\eprop

\begin{proof} By renumbering twigs we may assume $\bar u_1,\bar u_2,\bar u_3\leq\bar u_4,\bar u_5$, $d(T_1)\leq d(T_2)\leq d(T_3)$ and $\bar \delta_4\geq \bar \delta_5$. Suppose $\bar u_4, \bar u_5\neq 0$. We have $\bar \delta_4+\bar \delta_5\geq \frac{7}{2}-3(\frac{1}{2}+\frac{1}{3})=1$. Then $\bar \delta_4\geq \frac{1}{2}$, so $d(T_4)\leq 3$.

Suppose $d(T_4)=2$. By \ref{cor: chains with small u} $Q_4=[2,1,3,2]$, so $\frac{2}{d(T_5)}\geq \bar \delta_5\geq 1-\frac{7}{10}=\frac{3}{10}$ and $\bar u_5\leq \frac{1}{2}-\frac{1}{5}=\frac{1}{10}$. Thus $d(T_5)\leq 6$ and $\frac{1}{d(T_{10})}\geq \frac{3}{10}-\frac{1}{d(T_5)}$. By \ref{cor: chains with small u} we have $d(T_5)\in \{5,6\}$. Since $u([3,2]), u([2,3]), u([(2)_4]), u([(2)_5])\geq \frac{1}{5}$, by \ref{lem:chains with small d} we see that $T_5=[5]$ or $T_5=[6]$. If $T_5=[6]$ then we simultaneously have $d(T_{10})\leq 7$ and $T_{10}=[(2)_k,3,(2)_4]$ for some $k\geq 0$, which is impossible. Thus $T_5=[5]$ and then $d(T_{10})\leq 10$. This is possible only for $T_{10}=[2,2,2]$ or $T_{10}=[3,2,2,2]$. Then $\bar u_5=\frac{1}{2}$ and $\bar u_5=\frac{1}{3}$ respectively, a contradiction.

Suppose $d(T_4)=3$. If $T_4=[2,2]$ then by \ref{cor: chains with small u} $\bar \delta_5\geq 1-\frac{7}{12}=\frac{5}{12}$ (in particular $d(T_5)\leq 4$) and $\bar u_5\leq \frac{1}{2}-\frac{1}{3}=\frac{1}{6}$, which contradicts \ref{cor: chains with small u}. Thus $T_4=[3]$ and $Q_4=[3,1,2,3]$. We obtain $\bar \delta_5\geq \frac{7}{15}$ and $\bar u_5\leq \frac{3}{10}$, so $Q_5=[3,1,2,3]$ by \ref{cor: chains with small u}. By \eqref{eq:t=10, bounds on delta and e} $\bar \delta_1+\bar \delta_2+\bar \delta_3\geq \frac{7}{2}-2\cdot \frac{8}{15}=\frac{73}{30}$, so $\bar \delta_3\geq \frac{73}{30}-2\cdot \frac{5}{6}=\frac{23}{30}$, hence $d(T_1)=d(T_2)=d(T_3)=2$. On the other hand $\bar u_1+\bar u_2+\bar u_3\leq \frac{1}{2}-\bar u_4-\bar u_5=\frac{1}{10}$, so by \ref{cor: chains with small u} $\bar u_1=\bar u_2=\bar u_3=0$. This is the case (i).

We may now assume $\bar u_i=0$ for $i\neq 5$. Then \eqref{eq:t=10, bounds on delta and e} gives \begin{equation} \frac{1}{6}\leq \bar \delta_5\leq \bar e_5\leq \frac{2}{3},\label{eq:bounds on e5 delta5}\end{equation} where $\bar e_5=e(T_5)+e(T_{10})$. The lower bound $\bar \delta_5\geq \frac{1}{6}$ gives $d(T_5)\leq 11$. Since $d(T_5)$ and $d(T_{10})$ are coprime, we have in fact $\bar \delta_5>\frac{1}{6}$ and hence $\bar u_5<\frac{1}{2}$. In case $Q_5=[2,1,3]$ and in all cases listed in \ref{cor: chains with small u} the inequality $\bar e_5=\bar \delta_5+\bar u_5\leq\frac{2}{3}$ fails, hence $d(T_5)\geq 5$. It follows in particular that the twig $T_5$, cannot consists only of $(-2)$-curves, otherwise $\frac{2}{3}\geq e(T_5)=1-\frac{1}{d(T_5)}\geq \frac{4}{5}$.

Suppose $T_5=[n]$ for some $5\leq n\leq 11$. Then $T_{10}=[(2)_m,3,(2)_{n-2}]$ for some $m\geq 0$, so $d(T_{10})=(m+2)n-1$. The inequality $\bar u_5<\frac{1}{2}$ is equivalent to  $mn\leq 2$, hence $m=0$. For $n=10,11$ the inequality $\bar \delta_5>\frac{1}{6}$ fails, hence $n\leq 9$. This is the case (ii).

Since $d(T_5)\leq 11$, it is easy to list the remaining possibilities for $T_5$ with $e(T_5)<\frac{2}{3}$ (note that $e(T_5)\neq \frac{2}{3}$, as $e(T_{10})\neq 0$). These are: $[2,k]$ and $[k,2]$ for $3\leq k\leq 6$, $[3,3]$, $[3,4]$, $[4,3]$, $[2,3,2]$, $[3,2,2]$, $[4,2,2]$, $[2,3,2,2]$, $[3,2,2,2]$ and $[3,2,2,2,2]$. For each such $T_5$ we computed explicitly $T_{10}$ for which $\bar u_5<\frac{1}{2}$ and we checked that in cases other than $T_5=[5,2]$, $T_{10}=[3,2,2,3]$ one of the inequalities $\bar \delta_5>\frac{1}{6}$ or $\bar e_5\leq \frac{2}{3}$ fails.  \end{proof}

Decompose $\pi\:X\to \PP^2$ into blow-ups $\sigma_1\circ\ldots\circ\sigma_s$. Let $m_i$, $i=1,\ldots,s$ be the multiplicity if the center of $\sigma_i$ as a point on the respective proper transform of $\bar E$. The non-increasing sequence of multiplicities of a given cusp $q_i$ and of all centers lying above it is the \emph{multiplicity sequence of the cusp $q_i$}. Since $p_a(E)=0$, the genus formula reads as $$\sum_{i=1}^s\binom{m_i}{2}=\binom{\deg \bar E-1}{2}.$$

The multiplicity sequence for a cusp with $Q_i=[2,1,3]$ is $(2,1,1)$. For $Q_i=[3,1,2,3]$ it is $(3,2,1,1)$, for $Q_i=[n,1,(2)_{n-2},3]$ it is $(n,n-1,(1)_{n-1})$ and for $Q_i=[5,2,1,3,2,2,3]$ it is $(9,7,(2)_3,(1)_2)$. Thus $\binom{\deg \bar E-1}{2}=11$ in case (i), $\binom{\deg \bar E-1}{2}=n^2-2n+5$ in case (ii) and $\binom{\deg \bar E-1}{2}=64$ in case (iii). These equations have no solutions in natural numbers, which completes the proof of the theorem.

\section{Four cusps}\label{sec:four cusps}

Let $\bar E\subseteq \PP^2$ be a rational cuspidal curve with cusps $q_1,\ldots,q_c$. For $i=1\ldots,c$ let $\bar m_i=(m_{i,0},m_{i,1},\ldots,m_{i,k_i})$ be the multiplicity sequence of $q_i$ as defined above. Note that we do not omit $1$'s from the sequence. The resolution $\pi\:X\to \PP^2$ can be described in terms of Hamburger-Noether pairs (characteristic pairs). For a cusp $q_i$ we denote the sequence of H-N pairs by $$\binom{c_{i,1}}{p_{i,1}},\ldots,\binom{c_{i,h_i}}{p_{i,h_i}},$$ where $gcd(c_{i,j},p_{i,j})=c_{i,j+1}$ for $j<h_i$ and $gcd(c_{i,h_i},p_{i,h_i})=1$. In general the pairs depend on a choice of the system of local parameters $\{x_1,y_1\}$ around $q_i$. We choose our parameters so that $L_1=\{x_1=0\}$ is the tangent direction of $\bar E$. Recall that the first H-N pair is defined as $$c_{i,1}=(\bar E\cdot \{x_1=0\})_{q_i},\  p_{i,1}=(\bar E\cdot \{y_1=0\})_{q_i},$$ where $(\ )_{q_i}$ denotes the local intersection index at $q_i$. The inductive step defining $\{x_2,y_2\}$, and hence the remaining part of the sequence of H-N pairs for $q_i$, is as follows. We blow up over $q_i$ until the proper transform $E'$ of $\bar E$ meets the inverse image not in a node. Denote the last produced exceptional curve by $L_2$ and the point of intersection with $E'$ by $\wt q_i$. If $E'$ is smooth at $\wt q_i$ we put $h_i=1$. Otherwise we choose local parameters $\{x_2,y_2\}$ around $\wt q_i$, so that $\{x_2=0\}=L_2$ and $p_{i,2}=(E'\cdot \{y_2=0\})_{\wt q_i}$ is the multiplicity of $\wt q_i\in E'$. This forces $p_{i,j}\leq c_{i,j}$. It follows from the definition that $p_{i,1}=m_{i,0}$.

A more complete reference to H-N pairs is \cite{Russell2}. Here we basically followed \cite[1.12]{CKR-embeddings}.  By $\binom{c}{p}_k$ we mean a sequence of pairs $\binom{c}{p},\ldots,\binom{c}{p}$ of length $k$. For a short proof of the following result see for example \cite[2.10]{PaKo-general_type}.

\blem\label{lem:sum of m_i, m_i^2 using pairs} With the notation as above:\benum[(i)]

\item $\ds \sum_{j=1}^{k_i} m_{i,j}=c_{i,1}+\sum_{j=1}^{h_i}p_{i,j}-1,$

\item $\ds \sum_{j=1}^{k_i} m_{i,j}^2=\sum_{j=1}^{h_i}c_{i,j}p_{i,j}.$

\eenum \elem

Put $M(q_i)=\ds c_{i,1}+\sum_{j=1}^{h_i}p_{i,j}-1$ and $I(q_i)=\ds \sum_{j=1}^{h_i}c_{i,j}p_{i,j}$.

\bex If $Q_i=[2,1,3]$ then $h_i=1$, $\binom{c_{i,1}}{p_{i,1}}=\binom{3}{2}$, $M(q_i)=4$ and $I(q_i)=6$. If $Q_i=[5,2,1,3,2,2,3]$ then $h_i=1$, $\binom{c_{i,1}}{p_{i,1}}=\binom{16}{9}$, $M(q_i)=24$ and $I(q_i)=144$.  \eex

\bcor\label{cor:I and II} Let $\gamma=-E^2$ and let $d=\deg \bar E$. Then \benum[(i)]

\item $\gamma-2+3d=\ds \sum_i M(q_i),$

\item $\gamma+d^2=\ds \sum_i I(q_i).$ \eenum\ecor

\begin{proof} Let $C\subseteq X$ be an irreducible curve on a smooth projective surface. Let $p\in C$ be a singular point of $C$ having multiplicity $m$ and let $\sigma\: X'\to X$ be a blow-up at $p$. Denote the exceptional curve by $L$ and the proper transform of $C$ on $X'$ by $C'$. Then $K_X\cdot C=\sigma^*K_X\cdot C'=(K_{X'}-L)\cdot C'=K_X'\cdot C'-m$. Also $C^2=\sigma^*C\cdot C'=(C'+mL)\cdot C'=C'^2+m^2$. It follows that the sum of all multiplicities $m_{i,j}$ equals $K_X\cdot E-K_{\PP^2}\cdot \bar E=\gamma-2+3d$ and the sum of their squares equals $\bar E^2-E^2=d^2+\gamma$. \end{proof}

\bprop\label{prop:four cusps} If $\bar E\subseteq \PP^2$ is non-rectifiable and $c=4$ then $t\leq 9$. \eprop

\begin{proof} Suppose $t\geq 10$. By \ref{cor:nonrectifiable implies t<=10} $t=10$. Since $\kappa(2K+E)>0$ and $\kappa(K+2)=2$, \eqref{eq:*} and \eqref{eq:diamond} give $$\delta(D)=e(D)=4,$$ so all maximal twigs of $D$ are tips. In fact \eqref{eq:diamond} gives also $\cal P^2=3$ and $K\cdot (K+D)=1$. By \ref{lem:properties of u}(ii), renumbering the twigs if necessary, we may assume that $T_1$, $T_2$, $T_3$, $T_4$ are equal $[2]$ and $T_5$, $T_6$, $T_7$, $T_8$ are equal $[3]$ and that $T_i,T_{i+4}\subseteq Q_i$ for $i\leq 4$. Let $T_9\subseteq Q_1$ and $T_{10}\subseteq Q_1\cup Q_2$ be the remaining two maximal twigs of $D$. We may assume that in case $T_{10}\subseteq Q_2$ we have $d(T_9)\leq d(T_{10})$ and in case $T_{10}\subseteq Q_1$ the tip $T_9$ is created before $T_{10}$ in the process of resolving the singularity $q_1\in \bar E$. We have $\delta(T_{9})+\delta(T_{10})=\frac{2}{3}$, which gives $(T_9,T_{10})=([2],[6])$ or $(T_9,T_{10})=([6],[2])$ or $(T_9,T_{10})=([3],[3])$. The Noether formula for $X$ reads as $K^2+\#D=10$. We have $1=K\cdot (K+D)=K\cdot (K+D-E)+\gamma-2$, and $K\cdot Q_i=0$, $\#Q_i=3$ for $i=3,4$, so we obtain $$\#Q_1+\#Q_2=K\cdot Q_1+K\cdot Q_2+\gamma.$$

We define a subsequence of the sequence of characteristic pairs to be of type $*{(n,k)}$ if it is equal $\binom{\alpha n}{\alpha n}_k,\binom{\alpha n}{\alpha n-\alpha}$ for some $n\geq 2$, $\alpha>0$ and $k\geq 0$. Let $p\in C$ be a point on a smooth curve and let $\wt C$ and $C'$ be the total and proper transforms of $C$ after performing blowups over $p$ according to a sequence of type $*{(n,k)}$. Then $$\wt C=[-C'^2+1,(2)_{k-1},3,(2)_{n-2},1,n]$$ if $k\neq 0$ and $\wt C=[-C'^2+2,(2)_{n-2},1,n]$ otherwise. In particular this produces a $(-n)$-tip. We have $K\cdot \wt C-K\cdot C'=n-1$ and $\#\wt C-1=k+n$. The subsequences producing the tips $T_{9}$ and $T_{10}$ are of type $*{(d(T_9),k)}$ and $*{(d(T_{10}),l)}$ for some $k,l\geq 0$. Therefore $\#Q_1+\#Q_2=6+k+l+d(T_9)+d(T_{10})$ and $K\cdot Q_1+K\cdot Q_2=d(T_9)+d(T_{10})-2$. The Noether formula reduces now to \begin{equation}l=\gamma-k-8\label{eq:Noether}.\end{equation} We have the following five possibilities for sequences of H-N pairs.

\begin{case} $T_9=[2]$, $T_{10}=[6]$, $T_{10}\subseteq Q_1$. \end{case}

The sequence of pairs for $q_2$ is $\binom{3}{2}$ and for $q_1$ it is $\binom{36}{24}, \binom{12}{12}_k, \binom{12}{6}, \binom{6}{6}_l, \binom{6}{5}.$ Then $M(q_1)=70+12k+6l$, $I(q_1)=966 + 144 k + 36 l$, $M(q_2)=4$, $I(q_2)=6$. Using equations \ref{cor:I and II} and \eqref{eq:Noether} we eliminate $\gamma$ and $l$ and we obtain the relation $d^2-21d=444+66k$. Then $(2d+1)^2=6\mod 11$, which is impossible.

\begin{case} $T_9=[6]$, $T_{10}=[2]$, $T_{10}\subseteq Q_1$. \end{case}

The sequence of pairs for $q_2$ is as above and for $q_1$ it is $\binom{36}{24}, \binom{12}{12}_k, \binom{12}{10}, \binom{2}{2}_l, \binom{2}{1}.$ Then $M(q_1)=70+12k+2l$ and $I(q_1)=986 + 144 k + 4 l$. As above we get $d^2-9d=768+110k$, which gives $(2d+1)^2=3\mod 5$, a contradiction.

\begin{case} $T_9=[3]$, $T_{10}=[3]$, $T_{10}\subseteq Q_1$. \end{case}

The sequence of pairs for $q_2$ is as above and for $q_1$ it is $\binom{27}{18}, \binom{9}{9}_k, \binom{9}{6}, \binom{3}{3}_l, \binom{3}{2}$. Then $M(q_1)=52 + 9 k + 3 l$ and $I(q_1)=546 + 81 k + 9 l$. We get $d^2-12d=324+48k$, which gives $(d-6)^2=24\mod 48$, a contradiction.

\begin{case} $T_9=[3]$, $T_{10}=[3]$, $T_{10}\subseteq Q_2$. \end{case}

The sequences for $q_1$ and $q_2$ are $\binom{9}{6}, \binom{3}{3}_k, \binom{3}{2}$  and $\binom{9}{6}, \binom{3}{3}_l, \binom{3}{2}$ respectively. Then $M(q_1)=16 + 3 k$, $I(q_1)=60 + 9 k$, $M(q_2)=16+3l$ and $I(q_2)=60+9l$. We get $d^2-12d+12=0$, a contradiction.

\begin{case} $T_9=[2]$, $T_{10}=[6]$, $T_{10}\subseteq Q_2$.\end{case}

The sequences for $q_1$ and $q_2$ are $\binom{6}{4}, \binom{2}{2}_k, \binom{2}{1}$  and $\binom{18}{12}, \binom{6}{6}_l, \binom{6}{5}$ respectively. Then $M(q_1)=10 + 2 k$, $I(q_1)=26 + 4 k$, $M(q_2)=34 + 6 l$ and $I(q_2)=246 + 36 l$. We now eliminate $k$ and $l$ and we get $d^2-24d+52=-5\gamma$, hence $(2d+1)^2=3\mod 5$, a contradiction. \end{proof}

By \ref{thm:CN_holds_for_c>4}, \ref{cor:nonrectifiable implies t<=10} and \ref{prop:four cusps} we obtain the following result.

\bcor\label{cor:nonrectifiable implies t<=9} Let $\bar E\subseteq \PP^2$ be a rational cuspidal curve defined over complex numbers. If $\bar E\subseteq \PP^2$ is not rectifiable then the tree of the exceptional divisor for its minimal embedded resolution has at most nine maximal twigs. \ecor

\bibliographystyle{amsalpha}
\bibliography{bibl}

\providecommand{\bysame}{\leavevmode\hbox to3em{\hrulefill}\thinspace}
\providecommand{\MR}{\relax\ifhmode\unskip\space\fi MR }
\providecommand{\MRhref}[2]{%
  \href{http://www.ams.org/mathscinet-getitem?mr=#1}{#2}
}
\providecommand{\href}[2]{#2}
\begin{thebibliography}{CNKR09}

\bibitem[CNKR09]{CKR-embeddings}
Pierrette Cassou-Nogues, Mariusz Koras, and Peter Russell, \emph{Closed
  embeddings of {$\Bbb C\sp *$} in {$\Bbb C\sp 2$}. {I}}, J. Algebra
  \textbf{322} (2009), no.~9, 2950--3002.

\bibitem[Coo59]{Coolidge}
Julian~Lowell Coolidge, \emph{A treatise on algebraic plane curves}, Dover
  Publications Inc., New York, 1959.

\bibitem[Lan03]{Langer}
Adrian Langer, \emph{Logarithmic orbifold {E}uler numbers of surfaces with
  applications}, Proc. London Math. Soc. (3) \textbf{86} (2003), no.~2,
  358--396.

\bibitem[Laz04]{Lazarsfeld_II}
Robert Lazarsfeld, \emph{Positivity in algebraic geometry. {II}}, Ergebnisse
  der Mathematik und ihrer Grenzgebiete. 3. Folge. A Series of Modern Surveys
  in Mathematics [Results in Mathematics and Related Areas. 3rd Series. A
  Series of Modern Surveys in Mathematics], vol.~49, Springer-Verlag, Berlin,
  2004, Positivity for vector bundles, and multiplier ideals.

\bibitem[Miy01]{Miyan-OpenSurf}
Masayoshi Miyanishi, \emph{Open algebraic surfaces}, CRM Monograph Series,
  vol.~12, American Mathematical Society, Providence, RI, 2001.

\bibitem[MKM83]{Kumar-Murthy}
N.~Mohan~Kumar and M.~Pavaman Murthy, \emph{Curves with negative
  self-intersection on rational surfaces}, J. Math. Kyoto Univ. \textbf{22}
  (1982/83), no.~4, 767--777.

\bibitem[MT92]{MiTs-lines_on_qhp}
M.~Miyanishi and S.~Tsunoda, \emph{Absence of the affine lines on the homology
  planes of general type}, J. Math. Kyoto Univ. \textbf{32} (1992), no.~3,
  443--450.

\bibitem[Pal11]{Palka-exceptional}
Karol Palka, \emph{Exceptional singular {$\bold Q$}-homology planes}, Ann.
  Inst. Fourier (Grenoble) \textbf{61} (2011), no.~2, 745--774,
  \href{http://arxiv.org/abs/0909.0772}{arXiv:0909.0772}.

\bibitem[PK10]{PaKo-general_type}
Karol Palka and Mariusz Koras, \emph{Singular {$\bold Q$}-homology planes of
  negative {K}odaira dimension have smooth locus of non-general type}, to
  appear in Osaka J. Math.,
  \href{http://arxiv.org/abs/1001.2256}{arXiv:1001.2256}, 2010.

\bibitem[Rus80]{Russell2}
Peter Russell, \emph{Hamburger-{N}oether expansions and approximate roots of
  polynomials}, Manuscripta Math. \textbf{31} (1980), no.~1-3, 25--95.

\bibitem[Ton05]{Tono-on_the_number_of_cusps}
Keita Tono, \emph{On the number of the cusps of cuspidal plane curves}, Math.
  Nachr. \textbf{278} (2005), no.~1-2, 216--221.

\bibitem[Wak78]{Wakabayashi}
Isao Wakabayashi, \emph{On the logarithmic {K}odaira dimension of the
  complement of a curve in {$P\sp{2}$}}, Proc. Japan Acad. Ser. A Math. Sci.
  \textbf{54} (1978), no.~6, 157--162.

\end{thebibliography}
\end{document}